\documentclass[12pt]{article}

\usepackage{latexsym}
\usepackage{amsfonts}
\usepackage{amssymb}
\usepackage{amsmath}
\usepackage{mathrsfs}
\input xypic
\usepackage{hyperref}

\newcommand{\be}{\begin{equation}}
\newcommand{\ee}{\end{equation}}
\newcommand{\bea}{\begin{eqnarray}}
\newcommand{\eea}{\end{eqnarray}}
\newcommand{\bean}{\begin{eqnarray*}}
\newcommand{\eean}{\end{eqnarray*}}
\newcommand{\brray}{\begin{array}}
\newcommand{\erray}{\end{array}}

\newtheorem{dfn}{Definition}[section]
\newtheorem{thm}[dfn]{Theorem}
\newtheorem{lmma}[dfn]{Lemma}
\newtheorem{ppsn}[dfn]{Proposition}
\newtheorem{crlre}[dfn]{Corollary}
\newtheorem{xmpl}[dfn]{Example}
\newtheorem{rmrk}[dfn]{Remark}

\newcommand{\bdfn}{\begin{dfn}\rm}
\newcommand{\bthm}{\begin{thm}}
\newcommand{\blmma}{\begin{lmma}}
\newcommand{\bppsn}{\begin{ppsn}}
\newcommand{\bcrlre}{\begin{crlre}}
\newcommand{\bxmpl}{\begin{xmpl}}
\newcommand{\brmrk}{\begin{rmrk}\rm}

\newcommand{\edfn}{\end{dfn}}
\newcommand{\ethm}{\end{thm}}
\newcommand{\elmma}{\end{lmma}}
\newcommand{\eppsn}{\end{ppsn}}
\newcommand{\ecrlre}{\end{crlre}}
\newcommand{\exmpl}{\end{xmpl}}
\newcommand{\ermrk}{\end{rmrk}}

\author{ Jean Renault and S. Sundar }
\title{Groupoids associated to Ore semigroup actions. }

\begin{document}
\maketitle 
\begin{abstract}
 In this paper, we consider actions  of locally compact Ore semigroups on compact topological spaces. Under mild assumptions on the semigroup and the action, we construct a  semi-direct product groupoid with a Haar system. We also show that it is Morita-equivalent to a transformation groupoid. We apply this construction to the Wiener-Hopf $C^{*}$-algebras.
\end{abstract}
\footnote{\textbf{Acknowledgements}: 
This work was undertaken when the second author visited the first at Orl\'eans in October 2013. He thanks  Prof. Emmanuel Germain of the University of Caen and the GDR, 2947, NCG for the financial support which made the visit possible. He would also like to thank Anne Liger and the MAPMO (CNRS-UMR 7349) for  hospitality and CMI, Chennai for providing him leave.}

\noindent {\bf AMS Classification No. :} {Primary 22A22; Secondary 54H20, 43A65, 46L55.}  \\
{\textbf{Keywords.}} Wiener-Hopf C*-algebras, Semigroups, Groupoids.

\section{Introduction}

Because of developments in groupoid C*-algebras and a sustained interest for C*-algebras attached to semigroups and semigroup actions, it seems timely to revisit the groupoid approach to Wiener-Hopf C*-algebras initiated in \cite{Renault_Muhly}. Recall that the Wiener-Hopf operators on the real half-line are operators on the Hilbert space $L^2([0,\infty))$ given by
$$W(f)\xi(s)=\int_0^\infty f(s-t)\xi(t)dt$$
where $f\in L^1(\mathbb{R})$. The Wiener-Hopf C*-algebra ${\mathcal W}([0,\infty))$ of the real half-line is the C*-algebra generated by these operators. It has a particular simple structure which illuminates the basic index theory of these operators. Replacing the real half-line by a closed subsemigroup $P$ in a locally compact group $G$, one defines similarly the Wiener-Hopf C*-algebra ${\mathcal W}(G,P)$ (it is also called the Toeplitz algebra of $(G,P)$).  The article \cite{Renault_Muhly} studies the case when $P$ is a cone in $G= \mathbb{R}^n$, more specifically a polyhedral cone or a homogeneous self-dual cone. The groupoid description of these algebras gives their ideal structure; this was a first step towards index theorems for the corresponding Wiener-Hopf operators; the whole program, including index theorems and K-theory, has been successfully carried out recently by A. Alldridge and coauthor \cite{AJ07,AJ08,All11}.  Groupoids were used by A. Nica  in \cite{Nica_WienerHopf, Nica90, Nica92, Nica94} to study  the C*-algebra   of various subsemigroups of non-commutative goups. He also discovered that some conditions on the semigroup $P$ had to be imposed in order to obtain the groupoid description. In the discrete case, he introduced the class of quasi-lattice ordered semigroups, which includes the important example of the free semigroups. In the continuous case, he introduced a condition (M), which is satisfied by the closed convex cones with non-empty interior in  $\mathbb{R}^n$ and by the positive semigroup of the Heisenberg group. After establishing condition (M) for Lie subsemigroups of Lie groups, J. Hilgert and K.-H. Neeb  also applied in \cite{Hilgert_Neeb,Hilgert_Neeb1} groupoid techniques to study their Wiener-Hopf C*-algebras. 

In the present work, we show that a groupoid description is available for the Wiener-Hopf C*-algebras of an arbitrary locally compact Ore semigroup $P$. Our groupoid is constructed as the semidirect product $X\rtimes P$, where $X$ is the order compactification of $P$. The semidirect product construction was already used in Proposition 3.1 of \cite{Exel_Renault07} and requires the Ore condition. It generalizes a groupoid construction introduced in \cite{Renault_thesis} to describe the Cuntz algebra and developed by Deaconu (\cite{Dea95}). In \cite{Exel_Renault07}, the semigroup $P$ is assumed to be discrete; then the semidirect product is an \'etale groupoid. We show that, when the semigroup $P$ is locally compact, the semidirect groupoid is locally compact and has a continuous Haar system. Previous constructions introduced the Wiener-Hopf groupoid as a reduction of a semidirect product by a group. The present construction is more natural and avoids the delicate operation of groupoid reduction. However, it turns out that our groupoid is in fact a reduction of a semidirect product by a group (and agrees with the original one in the cases previously considered). We use first a general result, Theorem 1.8 of \cite{khoskam_Skandalis},  by Khoshkam and Skandalis to conclude that our groupoid is Morita equivalent to a semidirect product by a group. It seems appropriate to view this result as a dilation theorem when it is applied to a semidirect product by a semigroup. We also note that this is a particular case of the classical Mackey range construction. Because our semidirect product groupoid comes from an injective action, it can indeed be be viewed as a reduction of the dilation. These results  show that the Wiener-Hopf C*-algebras we consider are Morita equivalent to a reduced crossed product of a commutative C*-algebra by a group $G$, which makes the computation of its K-theory possible in some cases. 

Let us relate briefly the groupoid approach to Wiener-Hopf algebras to the general theory of semigroup $C^{*}$-algebras  and semigroup crossed products. The reduced C*-algebra of a discrete semigroup $P$ is usually defined as the C*-algebra generated by the Wiener-Hopf operators $\{ W(\delta_s): s \in P\}$; it agrees with our Wiener-Hopf algebra when $P$ is a discrete Ore semigroup. Several definitions of the universal C*-algebra of a discrete semigroup or semigroup action have been proposed. Let us quote some of the earliest articles on this subject (\cite{Peters84, Murphy87, Murphy91}) and refer the reader to \cite{DFK14} for a very complete recent survey. The universal $C^{*}$-algebra $C^{*}(G,P)$ of a quasi-lattice ordered  group $(G,P)$ introduced by Nica  in \cite{Nica92} is realised in \cite{LR96} as a crossed product by a semigroup of endomorphisms, however, contrarily to $\mathcal{W}(G,P)$, $C^{*}(G,P)$ is not given as a groupoid $C^{*}$-algebra (except for the amenable case where it agrees with $\mathcal{W}(G,P)$). In \cite{Exel_Renault07}, the full $C^{*}$-algebra $C^{*}(X  \rtimes P)$ of a semidirect product by a discrete Ore semigroup $P$ is shown to be an Exel crossed product.      Xin Li, motivated by the study of ring C*-algebras, gives in \cite{Li-semigroup}  a new definition of the full semigroup C*-algebra of a discrete cancellative semigroup (and of full crossed products of semigroup actions). Under weak assumptions which cover both Nica's quasi-lattice ordered semigroups and Ore semigroups, he shows in \cite{Li13} that the full semigroup C*-algebra admits a groupoid presentation.  A motivation for the present paper is to consider $C^{*}$-algebras associated to locally compact semigroups.

The organization of this paper is as follows.

In Section 2, we first introduce the class of locally compact semigroups considered in this work, namely the continuous Ore semigroups. Given such a semigroup $P$, its Wiener-Hopf C*-algebra is defined as the C*-algebra of operators on $L^2(P)$ generated by the Wiener-Hopf operators with symbol  in $L^{1}(P)$ (or in $L^{1}(G)$ if one prefers).

Section 3 gives the construction of the topological groupoid associated with the action of a topological  Ore semigroup on a compact topological space. This is the semi-direct product groupoid alluded to earlier. This work is limited to the case of an injective action. We give some basic examples.

In Section 4, our main Theorem 4.3 gives a very simple necessary and sufficient condition for the existence of a continuous Haar system for a semidirect product by a locally compact Ore semigroup. It combines techniques of \cite{Hilgert_Neeb} and ideas of \cite{Nica_WienerHopf}.

In Section 5, we show that the Wiener-Hopf C*-algebra of a continuous Ore semigroup $P$ is isomorphic to the reduced C*-algebra of the semi-direct product groupoid $X\rtimes P$, where $X$ is the order compactification of $P$. It also results from our study that continuous Ore semigroups satisfy Nica's condition (M). 

Section 6 contains the dilation theorem mentioned above and an alternative description of the Wiener-Hopf groupoid. For the sake of completeness, we recall the necessary background about groupoid equivalence and give the proof of Theorem 1.8 of \cite{khoskam_Skandalis} which characterizes the groupoids which are equivalent to semidirect products by groups. This is applied to the computation of the K-theory of some Wiener-Hopf C*-algebras.

Section 7 gives the example of the ``ax+b'' semigroup on $\mathbb{R}$.

All the topological spaces considered in this paper are assumed to be Hausdorff and second countable.

\section{The regular $C^{*}$-algebra of a continuous Ore semigroup}
Throughout this paper, we let $P$ to denote a closed subsemigroup of a locally compact group $G$ containing the identity element $e$. We also assume the following. 
\begin{enumerate}
 \item[(C1)] The group $G=PP^{-1}$, and
  \item[(C2)] The interior of $P$ in $G$, denoted $Int(P)$, is dense in $P$.
\end{enumerate}

\begin{rmrk}
Note that
\begin{itemize}
 \item (C1) is equivalent to the fact that $P$ generates $G$ and given $a,b \in P$, the intersection $aP \cap bP$ is non-empty. Such semigroups are called right reversible Ore semigroups.
\item Also note that $PInt(P)$ and $Int(P)P$ are contained in $Int(P)$. Thus $Int(P)$ is a semigroup. Since $Int(P)$ is nonempty, it follows  that $G=Int(P)Int(P)^{-1}$.
\end{itemize}
\end{rmrk}

Let $\mu$ be a left Haar measure on $G$ and let $\Delta$ be the modular function associated to $G$. Sometimes we write $dg$ for the Haar measure $\mu$. We always view $L^{2}(P)$ as a closed subspace of $L^{2}(G)$ by extending a function on $P$ to one on $G$ by declaring its value outside $P$ as zero. For $a \in P$, let $V_{a}$ be the operator on $L^{2}(P)$ defined as follows: For $\xi \in L^{2}(P)$, define $V_{a}(\xi)$ by the formula
\begin{displaymath}
V_{a}(\xi)(x)= \xi(xa)\Delta(a)^{\frac{1}{2}}.
\end{displaymath}
 For $\xi \in L^{2}(P)$, $V_{a}^{*}(\xi)$ is then given by 
\begin{displaymath}
\begin{array}{lll}
V_{a}^{*}(\xi)(x)&=& \left\{\begin{array}{lll}
                                 \xi(xa^{-1})\Delta(a)^{\frac{-1}{2}} &~if~x \in Pa \\
                                 0 & ~if~x \notin Pa .
                                 \end{array} \right. 
\end{array}
\end{displaymath}
For $g \in G$, let $U_{g}$ be the unitary on $L^{2}(G)$ defined by the right regular representation i.e. if $\eta \in L^{2}(G)$, then $U_{g}(\eta)(x)= \Delta(g)^{\frac{1}{2}}\eta(xg)$.

Let us denote the orthogonal projection of $L^{2}(G)$ onto $L^{2}(P)$ by $E$. Observe that for $a \in P$, $V_{a}=EU_{a}E$. Now the following are easily verifiable.

\begin{enumerate}
 \item[(1)] The maps $P \ni a \to V_{a} \in B(L^{2}(P))$ and $P \ni a \to V_{a}^{*} \in B(L^{2}(P))$ are strongly continuous,
 \item[(2)] For $a,b \in P$, $V_{a}V_{b}=V_{ab}$, 
 \item[(3)] For $a \in P$, $V_{a}^{*}$ is an isometry, and
 \item[(4)] If $g=ab^{-1}$, then $V_{a}V_{b}^{*}=EU_{g}E$.
\end{enumerate}

For $f \in L^{1}(P)$, let $W_{f}$ be the bounded operator on $L^{2}(P)$ given by \[W_{f}= \int_{a \in P}f(a)V_{a}~d\mu(a).\] The operator $W_{f}$ is the Wiener-Hopf operator associated to the symbol $f$. 

If $f \in L^{1}(G)$, we let \[W_{f}=\int_{g \in G}f(g)EU_{g}E~ d\mu(g).\] The $C^{*}$-algebra generated by $\{W_{f}: f \in L^{1}(G) \}$ is called the Wiener-Hopf algebra associated to $P$ and let us denote it by $\mathcal{W}(P)$.

\begin{ppsn}
\label{equality}
 Let $C_{red}^{*}(P)$ be the $C^{*}$-algebra generated by $\{W_{f}: f \in L^{1}(P)\}$. Then $C_{red}^{*}(P)=\mathcal{W}(P)$.
\end{ppsn}

We need the following lemma to prove Proposition \ref{equality}.

\begin{lmma}
\label{convolution}
 Let \[
      \mathcal{F}:= \{ \psi * \phi: \psi \in L^{1}(P),~\phi \in L^{1}(P^{-1})  \textrm{~and~} \psi, \phi \textrm{~have compact support} \}
     \]
Then the linear span of $\mathcal{F}$ is dense in $L^{1}(G)$.
\end{lmma}
\textit{Proof.} First note that if $K_{n}$ is a decreasing sequence of compact sets with non-empty interior such 
that $ \bigcap_{n}K_{n}=\{e\}$ then $ \displaystyle \chi_{n}:= \frac{1}{\mu(K_{n})}1_{K_{n}}$ is an approximate identity in $L^{1}(G)$.

Now let $U_{n}$ be a decreasing sequence of compact sets  such that  $e \in Int(U_{n})$ and 
$ \bigcap_{n}U_{n}=\{e \}$. Since $e \in Int(U_{n}) \cap P$ and $Int(P)$ is dense in $P$, it follows that $Int(U_{n}) \cap Int(P)$ is non-empty. Set $V_{n}:= U_{n} \cap P$. Then $V_{n}$ is a decreasing sequence of compact sets with non-empty interior such that $ \bigcap_{n}V_{n}=\{e\}$.
Let \[
     \phi_{n}:= \frac{1}{\mu(V_{n})}1_{V_{n}}.
    \]
Then $\phi_{n} \in L^{1}(P)$ and $(\phi_{n})$ is an approximate identity in $L^{1}(G)$. 

Consider  $f \in C_{c}(G)$  and let $K$ be its support. Since $G=(Int(P))(Int (P))^{-1}$, it follows that there exists $a_{1},a_{2}, \cdots ,a_{m} \in Int(P)$ such that \[\displaystyle K \subset \bigcup_{i=1}^{m}a_{i}Int(P)^{-1}.\]

By choosing a partition of unity, we can write $f=\sum_{i=1}^{m}f_{i}$ where $supp(f_{i}) \subset a_{i}(Int(P))^{-1}$. For each $i=1,2, \cdots m$, let \[\tilde{f_{i}}=L_{a_{i}^{-1}}f_{i}.\]where for $a \in G$ and $\phi \in C_{c}(G)$, we let $L_{a}(\phi)(x)=\phi(a^{-1}x)$. Then $\tilde{f_{i}} \in L^{1}(P^{-1})$. Also note that since $a_{i} \in P$, it follows that $L_{a_{i}}\phi_{n} \in L^{1}(P)$ for each $i$ and each $n$.

Observe that $ \mathcal{F} \ni \displaystyle \sum_{i=1}^{m}L_{a_{i}}\phi_{n}*\tilde{f_{i}} = \sum_{i=1}^{m}L_{a_{i}}(\phi_{n}*\tilde{f_{i}})$. Since left translations on $L^{1}(G)$ are bounded operators, it follows that as $n$ tends to infinity, $\sum_{i=1}^{m}L_{a_{i}}\phi_{n}* \tilde{f_{i}}$ converges to $\sum_{i=1}^{m}L_{a_{i}}\tilde{f_{i}}$ which equals $f$.

Thus we have proved that $C_{c}(G)$ is contained in the closure of the linear span of $\mathcal{F}$. Since $C_{c}(G)$ is dense in $L^{1}(G)$,  the proof is complete. \hfill $\Box$

Let $A:=\{(a,g) \in P \times G:g^{-1}a \in P\}$. Define $\sigma:A \to P \times P$ by $\sigma(a,g)=(a,g^{-1}a)$. Note that $\sigma$ is a bijection.
\begin{lmma}
 \label{Change of coordinates}
The Radon-Nikodym derivative of the push-forward measure $\sigma_{*}(\mu \times \mu|_{A})$ is given by 
\[
 \frac{d(\sigma_{*}(\mu \times \mu))}{d(\mu \times \mu)}(a,b)=\Delta(b)^{-1}.
\]
\end{lmma}
\textit{Proof.}
Let $\phi, \psi \in L^{1}(P)$ be positive. Now calculate to find that 
\begin{align*}
 \int \phi(a)\psi(g^{-1}a)&1_{P}(a)1_{P}(g^{-1}a)1_{A}(a,g)~ dadg \\
& =  \int \phi(a) 1_{P}(a) \Big ( \int \psi(g^{-1}a)1_{P}(g^{-1}a)1_{A}(a,g)~dg \Big) da \\
&= \int \phi(a) 1_{P}(a) \Big ( \int \psi(ga)1_{P}(ga)1_{A}(a,g^{-1})\Delta(g)^{-1} dg \Big) da \\
&= \int \phi(a) 1_{P}(a) \Big ( \int \psi(b)1_{P}(b)\Delta(b)^{-1}~db \Big) da \\
&=\int \phi(a)\psi(b) 1_{P}(a)1_{P}(b) \Delta(b)^{-1}~dadb .
\end{align*}
This completes the proof. \hfill $\Box$

\textit{Proof of Proposition \ref{equality}.} It is clear that $C_{red}^{*}(P) \subset \mathcal{W}(P)$. Let $\phi \in L^{1}(P)$ and $\psi \in L^{1}(P^{-1})$ be given. Define $\tilde{\phi},\tilde{\psi}$ as 
\begin{align*}
 \widetilde{\phi}(a):&=\phi(a)\\
 \widetilde{\psi}(a):&= \Delta(a)^{-1} \overline{\psi(a^{-1})}
\end{align*}
Then $\widetilde{\phi}, \widetilde{\psi} \in L^{1}(P)$.

Let $\xi, \eta \in L^{2}(P)$ be given. We have
\begin{align*}
 \langle W_{\widetilde{\phi}}W_{\widetilde{\psi}}^{*} \xi , \eta \rangle & =   \int \phi(a) \psi(b^{-1}) \langle V_{a}V_{b}^{*}\xi, \eta \rangle  \Delta(b)^{-1}~ da db \\
 & = \int \phi(a)\psi(a^{-1}g)1_{A}(a,g) \langle EU_{g}E \xi, \eta \rangle ~da dg\textrm{~~(by Lemma \ref{Change of coordinates})}\\
 &= \int \phi(a) \psi(a^{-1}g) \langle EU_{g}E \xi, \eta \rangle ~da dg \textrm{~~ (Since $\psi \in L^{1}(P^{-1})$)}\\
 & = \int (\phi * \psi )(g) \langle EU_{g}E \xi, \eta \rangle ~dg \\
 & = \langle W_{\phi * \psi}\xi, \eta \rangle 
\end{align*}
Thus it follows that $W_{\phi * \psi} \in C_{red}^{*}(P)$. Now by Lemma \ref{convolution}, the proof is complete.

\section{Semigroup actions and groupoids}
\begin{dfn}
 Let $X$ be a topological space and $P$ be a cancellative topological semigroup with identity $e$. A right action of $P$ on $X$ is a continuous map $X \times P \to X$, the image of $(x,a) \in X \times P$  is denoted by $xa$, such that 
\begin{enumerate}
 \item[(1)] For $x \in X$, $xe=x$, and
 \item[(2)] For $x \in X$ and $a,b \in P$, $(xa)b=x(ab)$.
\end{enumerate}
\end{dfn}

In what follows, we assume that $X$ is compact and $P$ satisfies the hypotheses of Section 2. Morover we assume that the action is injective i.e. for every $a \in P$, the map $X \ni x \to xa \in X$ is injective. This implies in particular that for $a \in P$, the map $X \ni x \to xa \in X$ is a homeomorphism from $X$ to $Xa$.
\begin{xmpl}
 Let $X=[0,\infty]$, the one point compactification of $[0,\infty)$ and $P=[0,\infty)$, the additive semigroup. The semigroup $P$ acts on the right on $X$ by translations. Here we understand that $\infty + a = \infty$.
\end{xmpl}
Before we discuss the next example, let us review the Vietoris topology.
 Let $(X,d)$ be a locally compact metric  space. Let $\mathcal{C}(X)$ be the collection of all closed subsets of $X$. Then $\mathcal{C}(X)$ endowed with Vietoris topology is  compact and is metrisable. We only need to know the convergence.  Let $(A_{n})$ be a sequence in $\mathcal{C}(X)$. Define 
\begin{align*}
 \lim \inf A_{n} & = \{ x \in X: \lim \sup d(x,A_{n}) = 0 \}, \textrm{~and~} \\
 \lim\sup A_{n} & = \{x \in X: \lim \inf d(x,A_{n})=0 \}.
\end{align*}

Then $A_{n}$ converges in $\mathcal{C}(X)$ if and only if $\lim\sup A_{n}=\lim \inf A_{n}$. Moreover if $\lim \sup A_{n}=\lim \inf A_{n}=A$ then $A_{n}$ converges to $A$. It is an easy exercise to show that if $F$ is a closed subset of $X$ then $\{A \in \mathcal{C}(X): A \subset F \}$ is closed. Or equivalently, if $U$ is open in $X$ then $\{A \in \mathcal{C}(X): A \cap U \neq \emptyset \}$ is open.

\begin{xmpl}
\label{Wiener_Hopf}
 Consider $\mathcal{C}(G)$, the space of closed subsets of $G$. We endow $\mathcal{C}(G)$ with the Vietoris topology. Then $G$ acts on $\mathcal{C}(G)$ on the right. For $A \in \mathcal{C}(G)$ and $g \in G$, $Ag = \{ag:a \in A\}$. Let $X$ be the closure of $\{P^{-1}a: a \in P \}$ in $\mathcal{C}(G)$. Then $P$ leaves $X$ invariant. The  space $X$ is called the order compactification of $P$.
\end{xmpl}

Let $X$ be a compact space on which $P$ acts on the right injectively. Let \begin{displaymath}
                                                                \mathcal{G}:= \{(x,g) \in X \times G: \exists ~a,b \in P,~ y \in X \textrm{~such that~}g=ab^{-1}, xa=yb \}.
                                                               \end{displaymath}
First let us prove that $\mathcal{G}$ is  closed in $X \times G$. We need a little lemma.
\begin{lmma}
\label{uniqueness}
 Let $g \in G$ and let $x,y \in X$. Suppose $g=ab^{-1}=a_{1}b_{1}^{-1}$. Then $xa=yb$ if and only if $xa_{1}=yb_{1}$.
\end{lmma}
\textit{Proof.}  Since $a_{1}b_{1}^{-1}=ab^{-1}$, it follows that $a_{1}^{-1}a=b_{1}^{-1}b$. Choose $\alpha, \beta \in P$ such that $a_{1}^{-1}a=\alpha \beta^{-1}$. Then $a_{1}\alpha=a \beta$ and $b_{1}\alpha = b \beta$. Now observe the following equivalences.
\begin{align*}
 xa_{1}=yb_{1} & \Leftrightarrow xa_{1}\alpha=yb_{1}\alpha ~~(\textrm{as the action  is injective})\\
               & \Leftrightarrow xa \beta = yb \beta \\
               & \Leftrightarrow xa=yb~~(\textrm{as the action  is injective}).
\end{align*}
This completes the proof. \hfill $\Box$

\begin{crlre}
\label{source}
 Let $g \in G$ and $x \in X$ be given. Then the following are equivalent.
\begin{enumerate}
 \item[(1)] $(x,g) \in \mathcal{G}$.
 \item[(2)] There exists a unique $y \in X$ such that if $g=ab^{-1}$ with $a,b \in P$ then $xa=yb$. We denote the unique element $y$ as $s(x,g)$.
\end{enumerate}
\end{crlre}
\textit{Proof.} Follows from Lemma \ref{uniqueness} and the injectivity of the action. \hfill $\Box$

Now we prove that $\mathcal{G}$ is closed. Let $(x_{n},g_{n})$ be a sequence in $\mathcal{G}$ such that $(x_{n},g_{n})$ converges to $(x,g)$. Write $g=ab^{-1}$. Choose open sets $U$, $V$ in $G$, having compact closures,  such that $(a,b) \in U \times V$. Then $UV^{-1}$ is an open set in $G$ containing $g$. Hence $g_{n} \in UV^{-1}$ eventually. Then $g_{n}=a_{n}b_{n}^{-1}$ with $a_{n} \in U$ and $b_{n} \in V$. Let $y_{n} \in X$ be such that $x_{n}a_{n}=y_{n}b_{n}$.  If necessary, by passing to a  subsequence of  $(a_{n},b_{n},y_{n})$, which is possible as $ U \times V \times X $ has compact closure in $G \times G \times X$, we can assume that $(a_{n},b_{n},y_{n})$ converges, say to $(a_{1},b_{1},y)$. Then the equality $x_{n}a_{n}=y_{n}b_{n}$ implies that $xa_{1}=yb_{1}$.

But $g_{n}=a_{n}b_{n}^{-1}$ converges to $g$. Hence $a_{1}b_{1}^{-1}=ab^{-1}$. Now by Lemma \ref{equality}, it follows that $xa=yb$. Thus $(x,g) \in \mathcal{G}$. Hence $\mathcal{G}$ is closed. 

We endow $\mathcal{G}$ with the subspace topology induced by the product topology on $X \times G$. Let $s:\mathcal{G} \to X$ be defined as follows. For $(x,g) \in \mathcal{G}$, $s(x,g)$ is the unique element in $X$ such that if $g=ab^{-1}$ with $a,b \in P$ then $xa=s(x,g)b$. The fact that $s$ is well defined follows from Corollary \ref{source}.

We claim that $s$ is continuous. Let $(x_{n},g_{n})$ be a sequence in $\mathcal{G}$ such that $(x_{n},g_{n})$ converges to $(x,g)$. Let $y_{n}=s(x_{n},g_{n})$. Since $X$ is compact, to show $s$ is continuous, it is enough to show that if $y_{n}$ converges to $y$ then $y=s(x,g)$. So suppose $y_{n}$ converges to $y$. As in the proof of the closedness of $\mathcal{G}$, we can write $g_{n}=a_{n}b_{n}^{-1}$ with $a_{n},b_{n}$ lying in a compact subset of $P$. Then $x_{n}a_{n}=y_{n}b_{n}$. Let $(a_{n_{k}},b_{n_{k}})$ be a convergent subsequence and let $(a,b)$ be its limit. Then it follows that $xa=yb$. But  $a_{n}b_{n}^{-1}=g_{n}$ converges to $g$. Therefore $g=ab^{-1}$. Hence $y=s(x,g)$. This proves that $s:\mathcal{G} \to X$ is continuous.

The space $\mathcal{G}$ has a groupoid structure. The multiplication and the inversion are given as follows.
 \begin{align*}
  (x,g)(y,h)&=(x,gh)~\textrm{if~} y=s(x,g)\\
  (x,g)^{-1}&=(s(x,g),g^{-1}).
 \end{align*}
Let us verify if $(x,g)$ and $(y,h)$ in $\mathcal{G}$ are composable, then $(x,gh) \in \mathcal{G}$. Let $g=ab^{-1}$, $h=cd^{-1}$ ( where $a,b,c,d \in P$) and $z=s(y,h)$. Choose $\alpha, \beta \in P$ such that $b^{-1}c=\alpha \beta^{-1}$. Then $gh=(a\alpha)(d\beta)^{-1}$.  Now
 \begin{align*}
  x(a\alpha) & = (xa)\alpha \\
             & = (yb)\alpha \\
             &= y(c \beta) \\
             & =(zd)\beta \\
             &= z(d\beta).
 \end{align*}
Hence $(x,gh) \in \mathcal{G}$. The inverse is well defined follows from the definition.

Since the map $s$ is continuous, it follows that the multiplication and the inversion are continuous. Since $\mathcal{G}$ is closed, $\mathcal{G}$ is locally compact. Thus $\mathcal{G}$ is a locally compact Hausdorff topological groupoid. The unit space  $\mathcal{G}^{(0)}$ is homeomorphic to $X$. We denote this groupoid $\mathcal{G}$ by $X \rtimes P$  and call $X \rtimes P$ as the semi-direct product of $X$ by $P$.
\section{Existence of Haar system}
The main aim of this section is to show that the topological groupoid $X \rtimes P$ has a left Haar system under some assumptions on the action. The proofs rely heavily on the techniques used in \cite{Hilgert_Neeb}.  We start with a lemma which is crucial to what follows. It is really Lemma II.12, page 97 in \cite{Hilgert_Neeb}. We include the proof for completeness.
\begin{lmma}
\label{Density}
 Let $P$ be a closed subsemigroup of a locally compact group $G$. Assume that $e \in P$ and $Int(P)$ is dense in $P$. Let $A \subset G$ be closed and $AP \subset A$. Then 
\begin{enumerate}
 \item[(1)] The interior of $A$, $Int(A)$, is dense in $A$, and
 \item[(2)] The boundary of $A$, $\partial A$, has measure $0$.
\end{enumerate}
\end{lmma}
\textit{Proof.}  Observe that $AInt(P)$ is an open set contained in  $A$. Hence $AInt(P) \subset Int(A)$. Since $Int(P)$ is dense in $P$, it follows that $AInt(P)$ and $Int(A)$ are dense in $A$.

Let $U$ be a open set containing  $e$ such that $U$ has compact closure.
 We claim that there exists a sequence $s_{n} \in U \cap Int(P)$ such that  $s_{n} \in Int(P)s_{n+1}$. 

  Since $U \cap P$ contains $e$, $U \cap Int(P)$ is non-empty. Choose $s_{1} \in U \cap Int(P)$. Assume that we have choosen $s_{1},s_{2},\cdots, s_{n}$. Now $U \cap Int(P) \cap P^{-1}s_{n}$ contains $s_{n}$. Since $(Int(P))^{-1}s_{n}$ is dense in $P^{-1}s_{n}$, it follows that $U \cap Int(P) \cap Int(P)^{-1}s_{n}$ is non-empty. Choose $s_{n+1} \in U \cap Int(P) \cap Int(P)^{-1}s_{n}$. This proves the claim.

To show that $\partial A$ has measure zero, it suffices to show that if $K \subset G$ is compact, then $\partial A \cap K$ has measure zero. Fix a right Haar measure $\lambda$ on $G$. 

Let $s_{n}$ be a sequence in $U \cap Int(P)$ such that $s_{n} \in Int(P)s_{n+1}$. Since $Int(P)$ is a semigroup, it follows that if $m>n$ then $s_{n} \in Int(P)s_{m}$. Since $AInt(P) \subset Int(A)$, it follows that if $m>n$, $\partial A s_{n} \subset Int(A)s_{m}$. Hence if $m>n$, $(\partial A \cap K)s_{n} \cap (\partial A \cap K)s_{m}$ is empty.

Now calculate to find that
 \begin{align*}
  \sum_{n=1}^{\infty}\lambda(\partial A \cap K)&= \sum_{n=1}^{\infty}\lambda((\partial A \cap K)s_{n}) \\
                                              &= \lambda (\bigcup_{n=1}^{\infty}(\partial A \cap K)s_{n}) \\
                                              & \leq \lambda(KU) \\
                                              & < \infty.
 \end{align*}

Hence $\lambda(\partial A \cap K)=0$. This completes the proof. \hfill $\Box$

For $x \in X$, let $\mathcal{G}^{x}:=r^{-1}(x)$. Here $r:\mathcal{G} \to X$ is the range map. For $x \in X$, let \[
                                                                                                                  Q_{x}:= \{g\in G: \exists ~a,b \in P \textrm{~such that~} g=ab^{-1} \textrm{~and~}  xa \in Xb\}.
                                                                                                                 \]
Note that $\mathcal{G}^{x}=\{x\} \times Q_{x}$ for $x \in X$. Also $g \in Q_{x}$ if and only if $(x,g) \in \mathcal{G}$. 

Since $\mathcal{G}$ is closed in $X \times G$, it follows that $Q_{x}$, being the preimage of $\mathcal{G}$ under the map $g \to (x,g)$, is a closed subset of $G$. Also note that for $z \in X$ and $a \in P$, $(z,a) \in \mathcal{G}$.
Now if $g \in Q_{x}$ and $a \in P$ then $(x,g)$ and $(s(x,g),a)$ are composable and hence its product $(x,ga) \in \mathcal{G}$. Thus $Q_{x}P \subset Q_{x}$. Now the following is an immediate corollary of Lemma \ref{Density}.

\begin{crlre}
\label{boundary}
 For $x \in X$, let $Q_{x}:= \{g \in G: (x,g) \in \mathcal{G}\}$. Then $Q_{x}$ is closed, its interior is dense in $Q_{x}$ and its boundary has measure $0$.
\end{crlre}

For $x \in X$, let $\lambda^{x}$ be the measure on $\mathcal{G}$ defined as follows. For $f \in C_{c}(\mathcal{G})$,
\[
 \int f d\lambda^{x} = \int f(x,g)1_{Q_{x}}(g)dg.
\]
Here $dg$ denotes the left Haar measure on $G$.

We can now state our main theorem.

\begin{thm}
\label{main theorem}
 The following are equivalent.
\begin{enumerate}
 \item[(1)] The groupoid $\mathcal{G}$ has a left Haar system.
 \item[(2)] The map $X \times Int(P) \to X$ is open.
 \item[(3)] If $U \subset Int(P)$ is open, then $XU$ is open in $X$.
 \item[(4)] For $a \in P$, $X Int(P)a$ is open in $X$.
 \item[(5)] The measures $(\lambda^{x})_{x \in X}$ form a left Haar system.
\end{enumerate}
\end{thm}

\textit{Proof.} Suppose $\mathcal{G}$ has a Haar system. Then the source  map $s:\mathcal{G} \to X$ is open. Observe that $X \times Int(P) \subset \mathcal{G}$. Also $X \times Int(P)$, being open in $X \times G$, is open in $\mathcal{G}$. Thus the action $X \times Int(P) \to X$, which is the source map, is open. This proves $(1)$ implies $(2)$.

The implication $(2) \Rightarrow (3)$ and $(3) \Rightarrow (4)$ are  obvious. Now let us prove $(4) \Rightarrow (5)$.

Assume $(4)$. For $x \in X$, $\overline{Int(Q_{x})}=Q_{x}$. Hence $supp(\lambda^{x})=\mathcal{G}^{x}$. Now we check that for $h \in C_{c}(\mathcal{G})$, the map $X \ni x \to \int h d\lambda^{x}$ is continuous. For $\phi \in C(X)$ and $f \in C_{c}(G)$, let $(\phi \otimes f)(x,g)=\phi(x)f(g)$.  Then $\phi \otimes f \in C_{c}(\mathcal{G})$ and the linear span of $\{ \phi \otimes f: \phi \in C(X), f \in C_{c}(G)\}$ is dense in $C_{c}(\mathcal{G})$ for the inductive limit topology. Hence it is enough to check that $X \ni x \to \int 1_{Q_{x}}(g)f(g)dg$ is continuous for every $f \in C_{c}(G)$.

Fix $f \in C_{c}(G)$. Let $(x_{n})$ be a sequence in $X$ such that $(x_{n})$ converges to $x$.

We claim that  $1_{Q_{x_{n}}} \to 1_{Q_{x}}$ a.e.

Suppose $g \notin Q_{x}$. Then $(x,g) \notin \mathcal{G}$. Since $\mathcal{G}$ is closed in $X \times G$, it follows that $(x_{n},g) \notin \mathcal{G}$ eventually. Hence  $g \notin Q_{x_{n}}$ eventually. Thus we have shown that $1_{Q_{x_{n}}}$ converges (pointwise) to $1_{Q_{x}}$ on  $G-Q_{x}$.

Now suppose $g \in Int(Q_{x})$. Let $g=ab^{-1}$ with $a,b \in P$. The map $G \times G \ni (\alpha,\beta) \to \alpha \beta^{-1} \in G$ is continuous. It follows that there exists open sets $U$,$V$ in $G$ such that $a \in U$, $b \in V$ and $UV^{-1} \subset Int(Q_{x})$. As $Int(P)b$ is dense in $Pb$ and $V \cap Pb$ contains $b$, it follows that the intersection $V \cap Int(P)b$ is non-empty. Let $c \in V \cap Int(P)b$. Then $ac^{-1} \in Q_{x}$. Thus $xa \in Xc$ which implies $xa \in X Int(P)b$. But $X Int(P)b$ is open in $X$ by assumption $(4)$. Also the sequence $x_{n}a \to xa$. Therefore $x_{n}a \in X Int(P)b$ eventually. Thus $x_{n}a \in Xb$ eventually which implies that eventually $(x_{n},g) \in \mathcal{G}$ or in other words $g \in Q_{x_{n}}$. This shows that $1_{Q_{x_{n}}}$ converges to $1_{Q_{x}}$ on $Int(Q_{x})$. 

Now by Corollary \ref{boundary}, $\partial Q_{x}$ has measure $0$. Thus $1_{Q_{x_{n}}} \to 1_{Q_{x}}$ a.e. 

By the dominated convergence theorem, it follows that $X \ni x \to \int 1_{Q_{x}}(g)f(g)dg$ is continuous for $f \in C_{c}(G)$.

Now let us verify left invariance. Let $(x,g) \in \mathcal{G}$ and denote $s(x,g)$ by $y$. Observe that $Q_{y}=g^{-1}Q_{x}$.  Consider a function $f \in C_{c}(\mathcal{G})$. Now 
 \begin{align*}
  \int f((x,g)(y,h))1_{Q_{y}}(h)dh & = \int f(x,gh)1_{g^{-1}Q_{x}}(h)dh \\
                                   & = \int f(x,u)1_{g^{-1}Q_{x}}(g^{-1}u)du \\
                                   & = \int f(x,u)1_{Q_{x}}(u)du.
 \end{align*}
This implies left invariance. Thus we have shown that $(4)$ implies $(5)$. The remaining implication $(5)$ implies $(1)$ is just definition. This completes the proof. \hfill $\Box$

 The above theorem should be compared to Proposition 1.3 of \cite{Nica_WienerHopf} which gives a necessary and sufficient condition such that the usual Haar system of a transformation groupoid $\mathcal {G}=U\rtimes G$ restricts to a Haar system of the reduced groupoid ${\cal G}_{V}$ where $V$ is a non-void locally closed subset of $U$. Under our assumptions, we obtain a much more convenient condition. 
\begin{crlre}
 Let $U \times G \to U$ be a continuous action where $U$ is a locally compact Hausdorff space. Let $V \subset U$ be a compact subset of $U$ and let $P:=\{g \in G: Vg \subset V\}$. Observe that $P$ is a closed subsemigroup of $G$. Suppose $G=PP^{-1}$ and $Int(P)$ is dense in $P$. Then the following are equivalent.
\begin{enumerate}
 \item[(1)] The usual Haar system on  $U \rtimes G$ reduces to a Haar system of the reduced groupoid  $(U \rtimes G)_{|V}$.
 \item [(2)] The map $V \times Int(P) \to V$ is open.
\end{enumerate}
\end{crlre}
\textit{Proof.} Observe that the reduction $(U \rtimes G)_{|V}=V \rtimes P$. Now the equivalence follows from Theorem \ref{main theorem}. \hfill $\Box$

We record the following useful fact from the proof of Theorem \ref{main theorem}.
\begin{rmrk}
\label{Uniqueness}
Let $X$ be a compact space on which $P$ acts on the right injectively. Suppose $\mathcal{G}:=X \rtimes P$ has a Haar system. Then the map $X \ni x \to Q_{x} \in \mathcal{C}(G)$ is continuous where $\mathcal{C}(G)$ is the space of closed subsets of $G$ with the Vietoris topology. For let $(x_{n})$ be a sequence in $X$ converging to $x$. Since $\mathcal{C}(G)$ is compact, we can assume that $Q_{x_{n}}$ converges and let $A$ be its limit. If $g \in A$ then there exists $g_{n} \in Q_{x_{n}}$ i.e. $(x_{n},g_{n}) \in \mathcal{G}$ such that $g_{n} \to g$. Since $\mathcal{G}$ is closed in $X \times G$ and $(x_{n},g_{n}) \to (x,g)$, it follows that $(x,g) \in \mathcal{G}$. Thus $g \in Q_{x}$ and $A \subset Q_{x}$. The proof of the implication $(4) \Rightarrow (5)$ of Theorem \ref{main theorem} implies that $1_{Q_{x_{n}}} \to 1_{Q_{x}}$ on $Int(Q_{x})$. Hence $Int(Q_{x}) \subset A$. But $A$ is closed and $Int(Q_{x})$ is dense in $Q_{x}$. Therefore we have $A=Q_{x}$. 
\end{rmrk}

The conditions in  Theorem \ref{main theorem} need not always have to be satisfied. Here is an example. 
Let $G:= \Big\{ \begin{bmatrix}
              a & b \\
              0 & 1
             \end{bmatrix} : a>0,~ b \in \mathbb{R} \Big \}.$ Then $G$ is isomorphic to the semidirect product $\mathbb{R} \rtimes \mathbb{R}^{*}_{+}$. Here the multiplicative group $\mathbb{R}^{*}_{+}$ acts on $\mathbb{R}$ by multiplication. Let $P:=[0,\infty) \times [1,\infty)$. Observe that the semigroup  $P$ is a closed in $G$, $PP^{-1}=G$ and $Int(P)$ is dense in $P$.

Let $\mathbb{CP}^{1}:=\mathbb{C} \cup \{\infty\}$ be the one-point compactification of $\mathbb{C}$. Consider the compact subset  $Y:=\{z \in \mathbb{CP}^{1}: Im(z) \geq 0\}$. Then $G$ acts on $Y$ on the right and the action is given by the formula \[
                                                                                                                            z.\begin{bmatrix}
                                                                                                                               a & b \\
0 & 1 \\
                                                                                                                              \end{bmatrix} = \frac{z-b}{a}.
\]
Here $\infty$ is left invariant. Let \[X:=\{z \in Y: Re(z) \leq 0, Im(z) \leq 1\} \cup \{\infty\}.\]    Then $P$ leaves $X$ invariant. The action of $P$ on $X$ is obviously injective. Note that $Int(P)=(0,\infty) \times (1,\infty)$. We leave it to the reader to verify that $X_{0}:=X(Int P)=\{(x,y): x<0, y<1\} \cup \{\infty\}$ which is not open in $X$. For $(-n,1) \notin X_{0}$ but converges to $\infty \in X_{0}$.                                                                                                           

\section{Wiener-Hopf $C^{*}$-algebras as groupoid $C^{*}$-algebras}

The main aim of this section is to show that Wiener-Hopf $C^{*}$-algebras can be realised as the reduced $C^{*}$-algebra of a groupoid. 

Let $X$ be a compact Hausdorff space on which $P$ acts on the right. Assume that the action $X \times Int(P) \to X$ is open so that the groupoid $X \rtimes P$ has a left Haar system. We endow $X \rtimes P$ with the left Haar system $(\lambda^{x})_{x\in X}$ as in Theorem \ref{main theorem}. Consider the following conditions.
\begin{enumerate}
 \item[(A1)] There exists $x_{0} \in X$ such that $Q_{x_{0}}=P$,
 \item[(A2)] The set $\{x_{0}a:a \in P\}$ is dense in $X$, and
 \item[(A3)] For $x,y \in X$, if $Q_{x}=Q_{y}$ then $x=y$.
\end{enumerate}
These conditions were introduced in a slightly different form by Nica in  \cite{Nica_WienerHopf} when he defines a generating system over $P$. The dilation of $X$ constructed in the next section is a minimal generating system over $P$ in his terminology. We show that if (A1), (A2) and (A3) are satisfied then the Wiener-Hopf $C^{*}$-algebra is isomorphic to $C_{red}^{*}(X \rtimes P)$. 

First we show that such a compact space exists. It turns out such a space is in fact unique. Recall the order compactification of $P$, from Example \ref{Wiener_Hopf}, which is the closure of $\{P^{-1}a: a \in P\}$ in the space of closed subsets of $G$ under the Vietoris topology. The semigroup $P$ acts on the order compactification by right multiplication. 
\begin{ppsn}
\label{Wiener-Hopf groupoid}
 Let $X$ be the order compactification of $P$. Then we have the following.
\begin{enumerate}
\item[(1)] The groupoid $X \rtimes P$ has a Haar system.
\item[(2)] For the pair $(X,P)$, the conditions (A1), (A2) and (A3) hold.
\item[(3)] Suppose $\widetilde{X}$ is a compact Hausdorff space with an injective right action of $P$. Suppose that $\widetilde{X} \rtimes P$ has a Haar system and (A1), (A2) and (A3) are satisfied then the map $\widetilde{X} \ni x \to Q_{x}^{-1} \in X$ is a $P$-equivariant homeomorphism.
\end{enumerate}
 \end{ppsn}
\textit{Proof.} To prove $(1)$, by Theorem \ref{main theorem}, it is enough to show that if $W \subset Int(P)$ is open then $XW$ is open in $X$. Let $W$ be an open subset of $Int(P)$. We claim that 
$XW=\{A \in X: A \cap W \neq \emptyset\}$. Suppose $A \in XW$. Then $A=Bw$ for some $w \in W$ and $B \in X$. Since $P^{-1} \subset B$, it follows that $w \in Bw=A$. Thus $A \cap W \neq \emptyset$.  Now suppose  $A \in X$ be such that  $A \cap W \neq \emptyset$.  Choose a sequence $a_{n} \in P$ such that $P^{-1}a_{n} \to A$. Choose $w \in A \cap W$. Then there exists $b_{n} \in P$ such that $b_{n}^{-1}a_{n} \to w$.  Since  $X$ is compact, by passing to a subsequence, if necessary, we can assume $P^{-1}b_{n}$ converges and let $B \in X$ be its limit. Since the action of $G$ on $\mathcal{C}(G)$ is continous, it follows that $(P^{-1}b_{n})b_{n}^{-1}a_{n}$ converges to $Bw$. But $P^{-1}b_{n}b_{n}^{-1}a_{n}=P^{-1}a_{n}$ converges to $A$. Hence $A=Bw$ for some $B \in X$ and $w \in W$. This proves the claim. By definition of the Vietoris topology, it follows that $\{A \in X: A \cap W \neq \emptyset\}$ is open in $X$. This proves $(1)$.

We claim that for $A \in X$, $Q_{A}=A^{-1}$. Let $A \in X$ be given. It is clear that $Q_{A}=\{g \in G: Ag \in X\}$.  Let $g \in Q_{A}$ be given. Then there exists a sequence $a_{n} \in P$ such that $P^{-1}a_{n} \to Ag$. Hence $P^{-1}a_{n}g^{-1} \to A$ in $\mathcal{C}(G)$. But then $a_{n}^{-1}a_{n}g^{-1} \to g^{-1}$. Thus $g^{-1} \in A$ or equivalently $g \in A^{-1}$. Now suppose $g \in A^{-1}$. Choose a sequence $a_{n} \in P$ such that $P^{-1}a_{n} \to A$. Then there exists a sequence $b_{n} \in P$ such that $b_{n}^{-1}a_{n} \to g^{-1}$. Since the action of $G$ on $\mathcal{C}(G)$ is continuous, it follows that $P^{-1}b_{n}=(P^{-1}a_{n})a_{n}^{-1}b_{n}$ converges to $Ag$. Hence $Ag \in X$ which implies $g \in Q_{A}$. This proves the claim.

Now take $x_{0}=P^{-1} \in X$. The fact that $Q_{A}=A^{-1}$ for $A \in X$ implies that (A1) and (A3) holds. Condition (A2) is just the definition of $X$. 

From Remark \ref{Uniqueness}, it follows that the map $\widetilde{X} \ni x \to Q_{x}^{-1} \in \mathcal{C}(G)$ is continuous.  Also observe that for $x \in X$ and $a \in P$, we have $Q_{xa}=a^{-1}Q_{x}$. This observation along with compactness of $\widetilde{X}$ and Conditions (A1) and (A2)   imply that the map $x \to Q_{x}^{-1}$ has range $X$ and is $P$-equivariant. Now (A3) gives the injectivity. Since $\widetilde{X}$ is compact, (3) follows. This completes the proof. \hfill $\Box$

We need the following proposition from \cite{Renault_Muhly}.

\begin{ppsn}[ Prop 2.17, \cite{Renault_Muhly}]
\label{Prop 2.17}
Let $\mathcal{G}$ be a groupoid with Haar system. Let $\mu$ be a measure on the unit space $\mathcal{G}^{(0)}$ such that $\mu(U)>0$ for every non-empty open invariant subset $U$ of $\mathcal{G}^{(0)}$. Then $Ind ~\mu$ is a faithful representation of $C_{red}^{*}(\mathcal{G})$.
\end{ppsn}

\begin{rmrk}
A part of the proposition quoted in  2.16 of \cite{Renault_Muhly} is wrong, namely the exactness of the sequence 
\[
0 \to  C_{r}^{*}(\mathcal{G}_{U}) \to C_{r}^{*}(\mathcal{G}) \to C_{r}^{*}(\mathcal{G}_{F}) \to 0
\]
where $U$ is an open invariant subset of $\mathcal{G}^{(0)}$ and $F=\mathcal{G}^{(0)}\backslash U$. It holds  for the full $C^{*}$-algebras. However the statement of Propositon \ref{Prop 2.17} is proved in complete detail in  2.17 of \cite{Renault_Muhly}.
\end{rmrk}

Now let $(X,P)$ be the unique dynamical system for which (A1), (A2) and (A3) holds. Consider the groupoid $\mathcal{G}:=X \rtimes P$ with the Haar system $(\lambda^{x})_{x \in X}$. For $f \in C_{c}(G)$, let $\widetilde{f}:\mathcal{G} \to \mathbb{C}$ be defined by $\tilde{f}(x,g)=f(g)$. Since $X$ is compact, it follows that $\widetilde{f} \in C_{c}(\mathcal{G})$. 

Consider the Dirac delta measure $\delta_{x_{0}}$ on $X$. Since the orbit of $\{x_{0}\}$ is dense, it follows,  from Proposition \ref{Prop 2.17},  that the induced representation $Ind(\delta_{x_{0}})$ gives a faithful representation of $C_{red}^{*}(\mathcal{G})$. Also observe that $\mathcal{G}^{x_{0}} = \{x_{0}\} \times P$ by (A1). Thus $L^{2}(\mathcal{G}^{x_{0}})$ is unitarily equivalent to $L^{2}(P)$. Now calculate to find that for $f \in C_{c}(G)$ and $\xi \in L^{2}(P)$,

 \begin{align*}
  Ind(\delta_{x_{0}})(\widetilde{f})(\xi)(a) &= \int \widetilde{f}((x_{0},a)^{-1}(x_{0},b))\xi(b) 1_{P}(b)db  \\
                                          &= \int f(a^{-1}b)\xi(b) 1_{P}(b) db \\
                                          & = \int f(u)\xi(au)1_{P}(au)du \\
                                          & = \int \Delta(u)^{\frac{-1}{2}}f(u) \Delta(u)^{\frac{1}{2}}\xi(au)1_{P}(au) du 
 \end{align*}
Thus it follows that $Ind(\delta_{x_{0}})(\widetilde{f})= W_{\widehat{f}}$ where $\widehat{f}(u)=\Delta(u)^{-\frac{1}{2}}f(u)$. Thus it follows that $\mathcal{W}(P) \subset C_{red}^{*}(\mathcal{G})$. To see the other inclusion we need the following lemma.

\textbf{Notations:} If $\phi \in C(X)$ and $f \in C_{c}(G)$,  let $\phi \otimes f \in C_{c}(\mathcal{G})$ be  defined by the equation $(\phi \otimes f)(x,g)=\phi(x)f(g)$. For $\phi \in C(X)$, $f \in C_{c}(G)$ and $a \in P$, let $R_{a}(\phi) \in C(X)$ be defined by $R_{a}(\phi)(x)=\phi(xa)$ and let $L_{a}(f) \in C_{c}(G)$ be defined by $L_{a}(f)(g)=f(a^{-1}g)$.

\begin{lmma}
\label{Density of Wiener-Hopf}
 Let $\mathcal{F}$ be a subset of $C(X)$ which is closed under conjugation. Assume that  $\mathcal{F}$  contains the constant function $1$ and the algebra generated by $\mathcal{F}$ is dense in $C(X)$. Then, with respect to the inductive limit topology on $C_{c}(\mathcal{G})$, the $*$-algebra generated by $\{ \phi \otimes f: \phi \in \mathcal{F}, f \in C_{c}(G)\}$ is dense in $C_{c}(\mathcal{G})$.
\end{lmma}
\textit{Proof.} Let $\phi_{1}, \phi_{2} \in C(X)$ and $f_{1},f_{2} \in C_{c}(G)$ be given. Now for $(x,g) \in \mathcal{G}$, 

\begin{align*}
\big((\phi_{1} \otimes f_{1})*(\phi_{2} \otimes f_{2})\big)(x,g) &= \phi_{1}(x) \int \phi_{2}(s(x,g).h)f_{1}(gh)f_{2}(h^{-1})1_{Q_{s(x,g)}}(h)dh \\
&=\phi_{1}(x) \int \phi_{2}(s(x,gh))f_{1}(gh)f_{2}(h^{-1})1_{Q_{x}}(gh)dh \\
& =\phi_{1}(x) \int \phi_{2}(s(x,u))f_{1}(u)f_{2}(u^{-1}g)1_{Q_{x}}(u)du.
\end{align*}
Let $a \in Int(P)$ be given. Choose a decreasing sequence of compact neighbourhoods $U_{n}$ containing $a$ in the interior such that $U_{n} \subset Int(P)$ and $ \bigcap_{n=1}^{\infty} U_{n}=\{a\}$. Let $f_{n} \in C_{c}(G)$  be such that $f_{n} \geq 0$, $supp(f_{n}) \subset U_{n}$ and $\int f_{n}(g) dg = 1$.
We claim that if $\phi,\psi \in C(X)$ and $f \in C_{c}(G)$ then  w.r.t  the inductive limit topology $(\phi \otimes f_{n})*(\psi \otimes f) \to \phi R_{a}(\psi) \otimes L_{a}(f)$. It is clear from   the above formula  that 
$(\phi \otimes f_{n}) * (\psi \otimes f)$ is supported in $X \times U_{1}supp(f)$ which is compact. 

Let $\epsilon >0$. Since $X$ is compact and $f$ is compactly supported, it follows that there exists $N \in \mathbb{N}$ such that for $n \geq N$, $(x,u,g) \in X \times U_{n} \times G$, we have
\[
 |\psi(xu)f(u^{-1}g)-\psi(xa)f(a^{-1}g)| \leq \epsilon.
\]
We leave the details of the proof to the reader. Now let $n \geq N$ be fixed. Now calculate as follows to find that for $(x,g) \in \mathcal{G}$
\begin{align*}
 |\big((\phi \otimes f_{n})*&(\psi \otimes f)\big)(x,g)-(\phi R_{a}\psi \otimes L_{a}(f))(x,g)| \\
 &= \big|\int_{u \in U_{n}}\phi(x)\psi(xu)f_{n}(u)f(u^{-1}g)du - \int_{u \in U_{n}}\phi(x)\psi(xa)f_{n}(u)f(a^{-1}g)du \big |\\
&\leq |\phi(x)|\int_{u \in U_{n}}|\psi(xu)f(u^{-1}g)-\psi(xa)f(a^{-1}g)|f_{n}(u)du \\
 &\leq \epsilon ||\phi||_{\infty}.
\end{align*}
This proves the claim.

Now let $\mathcal{A}$ be the $*$-algebra generated by $\{ \phi \otimes f:\phi \in \mathcal{F}, f \in C_{c}(G)\}$. Since $\mathcal{F}$ is dense in $C(X)$, it is enough to show that for $\phi_{1},\phi_{2}, \cdots ,\phi_{n} \in \mathcal{F}$ and $f \in C_{c}(\mathcal{G})$, $(\phi_{1} \phi_{2} \cdots \phi_{n}) \otimes f \in \overline{\mathcal{A}}$. We prove this by induction on $n$. For $n=1$, it is clear. 

Now let $\phi, \phi_{1},\phi_{2}, \cdots, \phi_{n} \in \mathcal{F}$ and $f \in C_{c}(G)$ be given. Set $\psi=\phi_{1}\phi_{2}\cdots \phi_{n}$. By induction hypothesis, it follows that $\psi \otimes f \in \overline{\mathcal{A}}$.

By our preceding discussion, it follows that for every $a \in Int(P)$, $\phi R_{a}\psi \otimes L_{a}f \in \overline{\mathcal{A}}$. Now let $a_{n}$ be a sequence in $Int(P)$ such that $a_{n} \to e$. Then $\phi R_{a_{n}}\psi \otimes L_{a_{n}}f \to \phi \psi \otimes f$. Thus it follows that $\phi \psi \otimes f \in \overline{\mathcal{A}}$.
This completes the proof. \hfill $\Box$

\begin{thm}
\label{maintheorem2}
 Let $X$ be a compact space on which $P$ acts on the right. Assume that the action $X \rtimes Int(P) \to X$ is open.
If $(X,P)$ satisfies (A1), (A2) and (A3) then $\mathcal{W}(P)$ is isomorphic to $C_{red}^{*}(X \rtimes P)$.
\end{thm}
\textit{Proof.} We have already shown that $\mathcal{W}(P) \subset C_{red}^{*}(X \rtimes P)$.  Clearly $1 \otimes f \in \mathcal{W}(P)$ for $f \in C_{c}(G)$. For $f \in C_{c}(G)$, let $\phi_{f} \in C(X)$ be given by $\phi_{f}(x)=\int 1_{Q_{x}}(t)f(t)dt$. Condition (A2) implies that the $\{\phi_{f}: f \in C_{c}(G)\}$ separates points of $x$.  This is because for every $x$, $Int(Q_{x})$ is dense in $Q_{x}$ and the boundary of $Q_{x}$ has measure zero. Thus if $1_{Q_{x}}=1_{Q_{y}}$ a.e. then $Q_{x}=Q_{y}$. Hence $\mathcal{F}:= \{\phi_{f}: f \in C_{c}(G)\}$ separates points of $X$ and the unital $*$-subalgebra generated by $\mathcal{F}$ is dense in $C(X)$. Now by Lemma \ref{Density of Wiener-Hopf}, to complete the proof, it is enough show that $\phi_{f} \otimes g \in \mathcal{W}(P)$ for every $f,g \in C_{c}(G)$. The proof of this is exactly the same as the proof of Proposition 3.5 in \cite{Renault_Muhly}.  One just have to replace $1_{X}(xt)$ in \cite{Renault_Muhly} by $1_{Q_{x}}(t)$.  Hence we omit the proof  \hfill $\Box$

We now indicate briefly, without proof, that the groupoid obtained in \cite{Renault_Muhly} and in Proposition \ref{Wiener-Hopf groupoid} are isomorphic. Let us recall the groupoid considered in \cite{Renault_Muhly}. Denote the algebra of uniformly continuous bounded functions on $G$ by $UC_{b}(G)$. Also $G$ acts on $UC_{b}(G)$ by right translation. Let $\widetilde{Y}$ be  the spectrum of the commutative $C^{*}$-subalgebra, denoted $\mathcal{A}$, of $UC_{b}(G)$ generated by $\{1_{P}*f: f \in C_{c}(G)\}$ where \[1_{P}*f(t)=\int 1_{P}(ts)f(s^{-1})ds= \int 1_{P^{-1}t}(s)f(s)\Delta(s)^{-1}ds.\]
Evaluation at points of $G$ gives multiplicative linear functionals of $\mathcal{A}$ and one obtains a $G$-equivariant map $\tau:G \to \widetilde{Y}$. Denote the closure $\tau(P)$ by $\widetilde{X}$ which is $P$-invariant. The groupoid considered in \cite{Renault_Muhly} is $\widetilde{X} \rtimes P$ or equivalently $\widetilde{Y} \rtimes G|_{\widetilde{X}}$.
For $a \in P$ and $f \in C_{c}(G)$, the equation
\[
 1_{P}*f(a)=\int 1_{P^{-1}a}(s)f(s)\Delta(s)^{-1}ds
\]
indicates that $\widetilde{X}$ is the closure of $\{1_{P^{-1}a}: a \in P\}$ in $L^{\infty}(G)$ where $L^{\infty}(G)$ is given the weak $^{*}$-topology. We show in the next proposition that the map $X \ni A \to 1_{A} \in L^{\infty}(G)$ is a $P$-equivariant embedding whose image is $\widetilde{X}$.

We finish this section by showing that Ore semigroups satisfies Condition (M) due to Nica. Let us recall the following definition due to Nica. Recall from \cite{Nica_WienerHopf} that a subset $A$ of $G$ is said to be solid if $Int(A)$ is dense in $A$ and the support of $\mu|_{A}$ is $A$. Here $\mu$ is a Haar measure on $G$.

\begin{dfn}[\cite{Nica_WienerHopf}]
A semigroup $P \subset G$ is said to satisfy Condition (M) if every element in the weak $^{*}$-closure of $\{1_{P^{-1}a}:a \in P\}$ in $L^{\infty}(G)$ is of the form $1_{A}$ for a solid closed subset $A$ of $G$.
\end{dfn}

\begin{ppsn}
\label{Condition M}
Let $P$ be a closed Ore semigroup of $G$ such that $Int(P)$ is dense in $P$. Denote the order compactification of $P$ by $X$.
\begin{enumerate}
\item[(1)] The map $X \ni A \to 1_{A} \in L^{\infty}(G)$ is a continuous $P$-equivariant embedding.
\item[(2)] The semigroup $P$ satisfies Condition (M).
\end{enumerate}
\end{ppsn}
\textit{Proof.} The continuity of the map $X \ni A \to 1_{A} \in L^{\infty}(G)$ follows from the fact  that $\{1_{Q_{A}}dg \}_{A \in X}$ is a Haar system for $X \rtimes P$ and $Q_{A}=A^{-1}$ for $A \in X$.  Observe that if $A \in X$ then $P^{-1}A \subset A$. Thus by Lemma \ref{Density}, it follows $Int(A)$ is dense in $A$ and the boundary of $A$ has measure $0$. Suppose $1_{A}=1_{B}$ a.e. for $A,B \in X$. Then $A \backslash B$ has measure zero. Hence the open set $Int(A)\backslash B$ has measure zero. This implies that $Int(A) \backslash B$ is empty. But $Int(A)$ is dense in $A$. Thus it follows that  $A\backslash B$ is empty.  Similarly $B \backslash A$ is empty.  This implies that the map $X \ni A \to 1_{A} \in L^{\infty}(G)$ is injective. We leave the $P$-equivariance to the reader. This proves $(1)$.

By (1), the weak$^{*}$-closure of $\{1_{P^{-1}a}: a \in P\}$ in $L^{\infty}(G)$ is $\{1_{A}: A \in X\}$. But for $A \in X$, $\overline{Int(A)}=A$. Thus $A$ is a solid closed subset for $A \in X$. This shows that $P$ satisfies Condition (M). This completes the proof. \hfill $\Box$

\section{Morita equivalence}

The groupoid ${\mathcal G}=X \rtimes P$ of the previous sections is defined as a semidirect product by the semigroup $P$. It admits a more usual presentation, namely as a reduction of a semidirect product $Y \rtimes G$ by the group $G$. As an intermediate step, we shall exhibit a Morita equivalence between $X \rtimes P$ and $Y \rtimes G$, where the $G$-space $Y$ is given by a classical construction, namely the Mackey range of a cocycle.

Let us first recall some definitions and notations concerning groupoid actions. A more detailed exposition can be found in \cite{MRW}. Let $\mathcal{G}$ be a locally compact groupoid with range and source maps $r,s:\mathcal{G} \to \mathcal{G}^{(0)}$, assumed to be open. Let $Z$ be a locally compact space on which $\mathcal{G}$ acts on the left. By definition, there is a map $\rho:Z \to \mathcal{G}^{(0)}$, assumed to be open and surjective and called  the moment map and an action map $(\gamma,z)\in\mathcal{G}*Z\rightarrow \gamma z\in Z$, where \[\mathcal{G}*Z:=\{(\gamma,z): s(\gamma)=\rho(z)\}\]
is the set of composable pairs.
The set $\mathcal{G}*Z$  has the structure of a groupoid, called the semidirect groupoid of the action and denoted by $\mathcal{G}\ltimes\, Z$, and  given by  
 \begin{align*}
  (\gamma,\gamma^{'}z)(\gamma^{'},z)&=(\gamma \gamma^{'},z)\\
   (\gamma,z)^{-1}&=(\gamma^{-1},\gamma.z) 
 \end{align*}
Endowed with the topology of the product $\mathcal{G}\times Z$, it becomes a topological groupoid. Let us check that the range and source maps are open, because it is the only fact we need. Since $s={\rm inv}\circ r$ and the inverse map ${\rm inv}((\gamma,z)=(\gamma^{-1},\gamma.z)$ is clearly a homeomorphism, it suffices to check that the range map is open. This map is simply the second projection $\pi_{2}:\mathcal{G}*Z \to Z$. If $U \times V$ is open in $\mathcal{G} \times Z$, then $\pi_{2}((U \times V)\cap (\mathcal{G}*Z))= \rho^{-1}(s(U)) \cap V$ is open.

Recall that the action is said to be free if $\gamma.z=z$ implies that $\gamma$ is a unit and proper if the map $\mathcal{G}*Z \ni (\gamma,z) \to (\gamma.z,z) \in Z \times Z$ is proper. If the action of $\mathcal{G}$ on $Z$ is proper, the quotient $\mathcal{G}\backslash Z$ is locally compact and Hausdorff.

Given two locally compact groupoids $\mathcal{G}$ and $\mathcal{H}$, a groupoid equivalence (also called Morita equivalence) is a locally compact space $Z$ which is a left $\mathcal{G}$-space, a right $\mathcal{H}$-space, the actions are free and proper, they commute and the corresponding moment maps $\rho$ and $\sigma$ identify respectively $\mathcal{G}\backslash Z\simeq \mathcal{H}^{(0)}$ and $Z/\mathcal{H}\simeq \mathcal{G}^{(0)}$.  Theorem 2.8 of \cite{MRW} says that, when ${\mathcal G}$ and ${\mathcal H}$ are endowed with Haar systems, $C_c(Z)$ can be completed into a $(C^*({\mathcal G}), C^*({\mathcal H}))$ imprimitivity bimodule. This is stated there for the full C*-algebras but the same result holds for the reduced C*-algebras.

A space endowed with a left $\mathcal{G}$-action and a right $\mathcal{H}$-action which commute is called a $(\mathcal{G},\mathcal{H})$-space. We shall encounter the following example of $(\mathcal{G},\mathcal{H})$-space. Let $j:\mathcal{G}\rightarrow\mathcal{H}$ be a continuous groupoid homomorphism. The space 
$$Z(j)=\{(x,\zeta)\in \mathcal{G}^{(0)}\times \mathcal{H}: j^{(0)}(x)=r(\zeta)\}$$
carries the left action $\gamma(s(\gamma),\zeta)=(r(\gamma), j(\gamma)\zeta)$ of $\mathcal{G}$ and the right action $(x,\zeta)\zeta'=(x,\zeta\zeta')$ of $\mathcal{H}$. Note that the right action of  $\mathcal{H}$ is free and proper. With an abuse of language, we call $Z(j)$ the graph of $j$.

Now we recall the Mackey range contruction. Let $\mathcal{G}$ be a locally compact groupoid with unit space $X=\mathcal{G}^{(0)}$, $G$ a locally compact group and $c:\mathcal{G}\to G$ a continuous cocycle (i.e. a continuous groupoid homomorphism). We consider its graph $Z=Z(c)=X\times G$ as above. The left action of $\mathcal{G}$ and the right action of $G$ are given by
$$\gamma(s(\gamma),a)=(r(\gamma), c(\gamma)a)\qquad (x,a)b=(x,ab).$$
As already mentioned, the right action of $G$ is free and proper, but not the left action of $\mathcal{G}$. Thus, in general the quotient $Y=\mathcal{G}\backslash Z$ is singular (for example, not Hausdorff). It is then natural to introduce the semidirect groupoid $\mathcal{G}\ltimes Z$, called the skew-product of the cocycle and denoted by $\mathcal{G}(c)$, as a substitute for this quotient space. The Mackey range of $c$ is defined in ergodic theory as the standard quotient $\underline Y=\mathcal{G}\backslash \backslash Z$ (i.e. the space of ergodic components), viewed as a right $G$-space. Under the assumption that the left action of $\mathcal{G}$ on $Z$ is proper, $Y=\mathcal{G}\backslash Z$ is a locally compact Hausdorff space endowed with a continuous action of $G$. Moreover, if the action of $\mathcal{G}$ on $Z$ is also free, then the skew-product $\mathcal{G}(c)$ is equivalent to the quotient space  $Y=\mathcal{G}\backslash Z$. We give now conditions under which the action of $\mathcal{G}$ on $Z$ is free and proper.  

\begin{dfn} \cite[Definition 1.6]{khoskam_Skandalis}. One says that a cocycle  $c:\mathcal{G}\to G$ is:
\begin{enumerate}
\item[(1)] faithful if the map from $\mathcal{G}$ to $\mathcal{G}^{(0)}\times G\times \mathcal{G}^{(0)}$ sending $\gamma$ to $(r(\gamma), c(\gamma), s(\gamma))$ is injective;
\item[(2)] closed if the above map from $\mathcal{G}$ to $\mathcal{G}^{(0)}\times G\times \mathcal{G}^{(0)}$ is closed;
\item[(3)]  injective if the map from $\mathcal{G}$ to $\mathcal{G}^{(0)}\times G$ sending $\gamma$ to $(r(\gamma), c(\gamma))$ is injective.
\end{enumerate}
\end{dfn}

One observes that, with above notation, $Z$ is a free $\mathcal{G}$-space if and only if $c$ is faithful and that $Z$ is a free and proper $\mathcal{G}$-space if and only if $c$ is faithful and closed. There is a slight abuse of language in (3); an equivalent definition is that $c^{-1}(e)\subset \mathcal{G}^{(0)}$. It is clear that a cocycle which is  injective is faithful and that the converse does not always hold.

\begin{thm}\cite[Theorem 1.8]{khoskam_Skandalis}  Let $\mathcal{G}$ be a locally compact groupoid with unit space $X=\mathcal{G}^{(0)}$, $G$ a locally compact group and $c:\mathcal{G}\to G$ a continuous cocycle. Assume that $c$ is faithful and closed. Then,
\begin{enumerate}
\item[(1)] the $(\mathcal{G},G)$-space $Z= Z(c)=X\times G$ is a groupoid equivalence between $\mathcal{G}$ and the semidirect product $Y\rtimes G$, where $Y=\mathcal{G}\backslash Z$ is the Mackey range of the cocycle;
\item[(2)] the map $j:\mathcal{G}\to Y\rtimes\,G$ such that  $j(\gamma)=([r(\gamma),e], c(\gamma))$, where $e$ is the unit element of $G$ and $[r(\gamma),e]$ is the class of $(r(\gamma),e)$ in $\mathcal{G}\backslash Z$, is a groupoid homomorphism; its graph $Z(j)$ is exactly $Z$ as a $(\mathcal{G}, Y\rtimes G)$-space;
\item[(3)] the image $X'=j^{(0)}(X)=\{[x,e], x\in \mathcal{G}^{(0)}\}$ of $\mathcal{G}^{(0)}$ in $Y$ meets every orbit under the action of $G$;
\item[(4)] if moreover $c$ is injective and $\mathcal{G}^{(0)}$ is compact, then $j$ is an isomorphism of $\mathcal{G}$ onto the reduction $(Y\rtimes\,G)_{|X'}$.
\end{enumerate} 
\end{thm}
\textit{Proof.} The right action of $G$ on $Z$ gives a right action of the semi-direct product $Y\rtimes G$: the moment map $\sigma: Z\rightarrow Y=\mathcal{G}\backslash Z$ is the quotient map. One defines  $z(\sigma(z), a)=za$ for $z\in Z$ and $a\in G$. This action remains free and proper. Since $c$ is faithful and closed, the left action of $\mathcal{G}$ is also free and proper. The identifications $X\simeq Z/Y\rtimes G$ and $Y\simeq \mathcal{G}\backslash Z$ are obvious. This proves (1). Let us check that the map $j$ of (2) is a groupoid homomorphism. Let $(\gamma,\gamma')\in \mathcal{G}^{(2)}$. We have
\begin{align*}
 j(\gamma)j(\gamma')&=([r(\gamma),e], c(\gamma))([r(\gamma'),e], c(\gamma'))\\
  &=([r(\gamma),e], c(\gamma))([\gamma(s(\gamma),e)], c(\gamma'))\\
&=([r(\gamma),e], c(\gamma))([r(\gamma),c(\gamma)], c(\gamma'))\\
&=([r(\gamma),e], c(\gamma))([r(\gamma),e]c(\gamma), c(\gamma'))\\
&=([r(\gamma\gamma'),e], c(\gamma\gamma'))=j(\gamma\gamma').
\end{align*}
The map $j^{(0)}: X\rightarrow Y=\mathcal{G}\backslash Z$ is given by $j^{(0)}(x)=[x,e]$. Therefore, we can identify
$Z(j)=\{(x,(y,a))\in X\times (Y\times G): j^{(0)}(x)=y\}$  and $Z=X\times G$ by sending $(x,([x,e],a))$ to $(x,a)$. Then, it is straightforward to check that the left actions of $\mathcal{G}$ [resp. right actions of $Y\rtimes G$] are the same. Since for all $(x,a)\in X\times G$, we have $(x,a)=(x,e)a$, the assertion (3) is true. Let us finally prove (4). By construction, $j$ is injective if and only if $c$ is injective.  Let us check that the image of $j$ is $(Y\rtimes G)_{|X'}$. The range of $j(\gamma)$ is $[r(\gamma),e]$ and its source is $[s(\gamma),e]$ which both belong to $X'$. Conversely, suppose that $(y,a)$ belongs to $(Y\rtimes G)_{|X'}$. Since $y\in X'$,  it is of the form $[x,e]$ with $x\in X$. Since $ya=[x,a]\in X'$, there exists $\gamma\in \mathcal{G}$ such that $(x,a)=\gamma(s(\gamma),e)$. This implies that $(y,a)=([r(\gamma),e], c(\gamma))=j(\gamma)$. If $X=\mathcal{G}^{(0)}$ is compact, $X'=j^{(0)}(X)$ is compact $j^{(0)}$ is a homeomorphism of $X$ onto $X'$. The compactness of $X$ and the closedness of $c$ imply that the map $(r,c): \mathcal{G}\rightarrow X\times G$ is closed. The map $j^{(0)}\times {\rm id}: X\times G\rightarrow Y\times G$ is also closed. Hence $j: \mathcal{G}\rightarrow Y\times G$ which is the composition of these maps is also closed. Since it is injective, it is a homeomorphism onto a closed subset of $Y\times G$. This completes the proof. \hfill $\Box$

\begin{rmrk}
Theorem 1.8 of \cite{khoskam_Skandalis} gives the assertions (1) and (2) (and not (3) and (4)) of the above theorem. We have given the proof in full (essentially the same as the original proof) to emphasize that its key point is the Mackey range construction, although it does not appear under that name in the original proof. Moreover the authors make a third assumption on the cocycle: they require it to be transverse. But this is exactly the condition that the source map of the skew-product $\mathcal{G}(c)$ is open. We have seen that this condition is automatically satisfied (provided that, as usual the range and source maps of $\mathcal{G}$ are open). 
\end{rmrk}

We apply this to the semi-direct product $\mathcal{G}=X\rtimes P$, where $X$ is a compact Hausdorff space on which $P$ acts on the right.  We assume that the action $X \rtimes Int(P) \to X$ 
is an open map or equivalently $X\rtimes P$ has a Haar system. Then the range and source maps of $X\rtimes P$ are open. The canonical cocycle $c: X\rtimes P\to G$, which is given by $c(x,ab^{-1})=ab^{-1}$ is faithful and closed. Therefore $\mathcal{G}$ is equivalent to the Mackey range semi-direct product $Y\rtimes G$, where $Y= \mathcal{G}\backslash X\times G$. In fact, since $c$ is injective and we assume that $X$ is compact, $\mathcal{G}$ is isomorphic to the reduction $(Y\rtimes\,G)_{|X'}$, where $X'=\{[x,e], x\in X\}$. 

\begin{ppsn}
\label{dilation} Under the above assumptions, let $X_{0}=X Int(P)$ and $j^{(0)}:X \to Y$ be the embedding given by the theorem.
\begin{enumerate}
 \item[(1)] The embedding $j^{(0)}:X \to Y$ is $P$-equivariant,
 \item[(2)] The image $j^{(0)}(X_{0})$ is open in Y, and
 \item[(3)] $Y=\bigcup_{a \in P}j(X)a^{-1}=\bigcup_{a \in Int(P)}j^{(0)}(X_{0})a^{-1}$.
 \item[(4)] Conditions $(1)$, $(2)$ and $(3)$ uniquely determines $Y$ up to a $G$-equivariant isomorphism. 
\end{enumerate}
\end{ppsn}

\textit{Proof.}
Observe that for $x \in X$ and $a \in P$, $(x,a) \in X \rtimes P$ and $(x,a).(xa,e)=(x,a)$ i.e. $[(xa,e)]=[(x,a)]$. Since $G=PP^{-1}$, it follows that $Y=\bigcup_{a \in P}X^{'}a^{-1}$. Now we verify that the image of $X_{0}:=X(Int(P))$ in $Y$ is open. It is equivalent to showing that \[A:= \{(x,ab^{-1}) \in X \times G: \textrm{~there exists~} y \in X_{0}, ~xa=yb\}\] is open in $X \times G$. Let $(x,g) \in A$ be given.  Write $g=ab^{-1}$ with $a,b \in Int(P)$. Let $y \in X_{0}$ be such that $xa=yb$. Then $xa \in X(Int P)b$ which is open in $X$. Thus there exist an open set $U \times V \subset X \times Int(P)$ containing $(x,a)$ and  $UV \subset X_{0}b$. Then $U \times Vb^{-1} \subset A$ and contains $(x,g)$. Thus $A$ is open in $X \rtimes G$. Since $G=Int(P)Int(P)^{-1}$, it is clear that $Y=\bigcup_{a \in Int(P)}j^{(0)}(X_{0})a^{-1}$. 
We leave (4) to the reader.
 \hfill $\Box$

The Morita equivalence established in this section  is useful to compute the $K$-groups. For example, we have the following version of Connes-Thom isomorphism for solid closed convex cones. Let $P \subset \mathbb{R}^{n}$ be a closed convex cone (i.e. $P$ is convex and $\alpha v \in P$ if $\alpha \geq 0$ and $v \in P$). Assume that $P$ spans $\mathbb{R}^{n}$. Then by duality theory, it follows that $Int(P)$ is dense in $P$.  We use additive notation for actions of $P$.

\begin{ppsn}
 \label{Connes-Thom}
Let $P$ be a  closed convex cone in $\mathbb{R}^{n}$ such that $P-P=\mathbb{R}^{n}$. Let $X$ be a compact Hausdorff space on which $P$ acts and assume that the action $X \times Int(P) \to X$ is open. Set $X_{0}=X+ Int(P)$. Then $K_{i}(C^{*}(X \rtimes P))$ is isomorphic to $K_{i+n}(C_{0}(X_{0}))$ for $i=0,1$.
\end{ppsn}
\textit{Proof.} Let $Y$ be the dilation of $X$ as in Proposition \ref{dilation}. For $s \in \mathbb{R}^{n}$, let $L_{s}$ be the translation on $C_{0}(Y)$  defined by $L_{s}(f)(y)=f(y-s)$ for $f \in C_{0}(Y)$ and $y \in Y$. We consider $C_{0}(X_{0})$ as a subset of $ C_{0}(Y)$. Then for $s \in P$, $L_{s}$ leaves $C_{0}(X_{0})$ invariant.

As  $X \rtimes P$ is Morita equivalent to $Y \rtimes \mathbb{R}^{n}$, it follows that thus $K_{i}(C^{*}(X \rtimes P))$ is isomorphic to $K_{i}(C_{0}(Y) \rtimes \mathbb{R}^{n})$. Now by Connes-Thom isomorphism, it follows that $K_{i}(C^{*}(X \rtimes P))$ is isomorphic to $K_{i+n}(C_{0}(Y))$. Thus it remains to show that for $i=0,1$, $K_{i}(C_{0}(Y)) \cong K_{i}(C_{0}(X_{0}))$. 

For $a \in Int(P)$, let $A_{a}$ be the closure of $\{ f \in C_{c}(Y): supp(f) \subset X_{0}-a\}$.  Since $(Int(P), <)$ is directed, where we write $a<b$ if $b-a \in Int(P)$, it follows that $\displaystyle \bigcup_{a \in Int(P)}A_{a}$ is dense in $C_{0}(Y)$. Also if $a<b$ then $A_{a} \subset A_{b}$. Thus $C_{0}(Y)$ is the inductive limit of $\displaystyle (A_{a})_{a \in Int(P)}$. Clearly $A_{a}$ is isomorphic to $ C_{0}(X_{0}-a) \cong C_{0}(X_{0})$. 

For $a,b \in Int(P)$ with $a<b$, let $i_{b,a}:A_{a} \to A_{b}$ be the inclusion. Under the isomorphism $A_{a} \cong C_{0}(X_{0})$, the map $i_{b,a}:C_{0}(X_{0}) \to C_{0}(X_{0})$ is nothing but $L_{b-a}$. Then $(L_{t(b-a)})_{t \in [0,1]}$ is a homotopy of homomorphisms on $C_{0}(X_{0})$ connecting the identity map to $L_{b-a}$. Thus the connecting map $i_{a,b}$ induce the identity map at the $K$-theory level. The proof is now complete by appealing to the continuity of $K$-groups under inductive limits. \hfill $\Box$

\begin{rmrk}
 Proposition \ref{Connes-Thom} gives a conceptual explanation of the vanishing of $K$-groups of the classical Wiener Hopf algebras associated to the additive semigroup $[0,\infty)$. In  this case $X=[0,\infty]$, the one point compactification of $[0,\infty)$, and $X_{0}=(0,\infty] \cong (0,1]$ whose $K$-groups are trivial. Similar observations have been made in \cite{KS97}. The index theorems associated to Wiener-Hopf operators (associated to polyhedral cones) have been studied extensively in \cite{AJ07}, \cite{AJ08} and \cite{All11}. In particular, it is established that the $K$-theory of the Wiener-Hopf algebra associated to polyhedral cones is trivial. 
\end{rmrk}

\section{An example}
We finish this article with an example. Let \[G:=\Big \{\begin{pmatrix}
                                                   a & b \\
                                                   0 & 1
                                                  \end{pmatrix}: a > 0 , b \in \mathbb{R} \Big\}.\] The group $G$ is isomorphic to the semidirect product $\mathbb{R} \rtimes (0,\infty)$. Let  $P= [0,\infty) \rtimes [1,\infty)$. 
Then $P$ is a semigroup and we leave it to the reader to verify that $PP^{-1}=P^{-1}P=G$. Also observe that $Int(P)=(0,\infty) \rtimes (1,\infty)$.

Let $Y:=[-\infty,\infty) \times [0,\infty]$. The group $G$ acts on $Y$ on the right and the action is given by the formula  , let \[
                                 (x,y).\begin{pmatrix}
        a & b \\
        0 & 1                                                                                                                
     \end{pmatrix}=\Big(\frac{x-b}{a},\frac{y}{a}\Big).
                                \]
Let $X:=[-\infty,0] \times [0,1]$. Then $X$ is $P$-invariant i.e. $XP \subset X$. Let $X_{0}:=X(Int P)$. We leave it to the reader to verify that  $X_{0}=[-\infty,0) \times [0,1)$.  Thus $X_{0}$ is open in $Y$ and hence in $X$. Thus for every $\gamma \in P$, $X_{0}\gamma$ is open in $Y$ and consequently in $X$. By Theorem \ref{main theorem}, it follows that the groupoid $X \rtimes P$ has a Haar system. We now verify the conditions (A1), (A2) and (A3) of Section 5.  

For $(x,y) \in X$, let $Q_{(x,y)}:=\{ g \in G: (x,y)g \in X\}$. Then, by definition, it follows that\[
                                                                            Q_{(x,y)}=\Big \{\begin{pmatrix}
                                                                                            a & b \\
                                                                                            0 & 1
                                                                                           \end{pmatrix} : x \leq b,~ y \leq a \Big \}.
                                                                           \]
Let $x_{0}=(0,1) \in X$. Then  $Q_{x_{0}}=P$. Note that the $P$-orbit of $x_{0}$ is $(-\infty,0] \times (0,1]$ which  is dense in $X$. Thus  (A1) and (A2) are satisfied.

Let $(x_{1},y_{1}),(x_{2},y_{2}) \in X$ be such that $Q_{(x_{1},y_{1})}=Q_{(x_{2},y_{2})}$. Suppose $x_{1}<x_{2}$. Choose $b \in \mathbb{R}$ be such that $x_{1}<b<x_{2}$. Then $\begin{pmatrix}
                                  1 & b \\
0 & 1
        \end{pmatrix}$ is in $Q_{(x_{1},y_{1})}$ but not in $Q_{(x_{2},y_{2})}$. Hence $x_{1} \geq x_{2}$. Similarly $x_{2} \geq x_{1}$. Thus $x_{1}=x_{2}$. 
        
Suppose $y_{1}<y_{2}$. Choose $a >0$ be such that $y_{1}<a<y_{2}$. Then $\begin{pmatrix}
                                                                          a & 0\\
                                                                          0 & 1
                                                                         \end{pmatrix}$ is in $Q_{(x_{1},y_{1})}$ but not in $Q_{(x_{2},y_{2})}$. Thus $y_{1} \geq y_{2}$. Similarly $y_{2} \geq y_{1}$. Hence $y_{1}=y_{2}$. Thus $(x_{1},y_{1})=(x_{2},y_{2})$. Hence (A3) is satisfied.

Also observe that $Y=\displaystyle \bigcup_{\gamma \in Int(P)}X_{0}\gamma^{-1}=\displaystyle \bigcup_{\gamma \in P}X \gamma^{-1}$. Thus by Theorem \ref{maintheorem2} and the Morita equivalence established in Section 6, it follows that the Wiener-Hopf $C^{*}$-algebra associated to $P$ is isomorphic to $C_{red}^{*}(X \rtimes P)$ and is Morita equivalent to the crossed product $C_{0}(Y) \rtimes (\mathbb{R} \rtimes (0,\infty))$.

So far we have only considered the right regular representation. One could also consider the left regular representation. Suppose that $P$ is a closed subsemigroup of $G$. Assume that $PP^{-1}=P^{-1}P=G$ and $Int(P)$ is dense in $P$. For $g \in G$, let $L_{g}$ be the unitary on $L^{2}(G)$  defined by $L_{g}(f)(x)=f(g^{-1}x)$. Let us denote the $C^{*}$-algebra on $L^{2}(P)$ generated by $\{ \int_{g \in P}f(g)EL_{g}E ~dg:~f \in C_{c}(G)\}$ by $\mathcal{W}_{\ell}(P)$. Here $E$ denotes the projection onto $L^{2}(P)$. We denote the Wiener-Hopf algebra associated to the right regular representation by $\mathcal{W}_{r}(P)$.  Consider  the unitary $U:L^{2}(G) \to L^{2}(G)$ defined by $(Uf)(x):=f(x^{-1})\Delta(x^{-1})^{\frac{1}{2}}$. The conjugation by $U$ maps $\mathcal{W}_{\ell}(P)$ onto $\mathcal{W}_{r}(P^{-1})$.

Let us now return to the example of $P=[0,\infty) \rtimes [1,\infty) \subset G =\mathbb{R} \rtimes (0,\infty)$. Observe that $P^{-1}=(-\infty,0] \rtimes (0,1]$. The map $G \in (b,a) \to (-b,a) \in G$ is an isomorphism, preserves the Haar measure and sends $P^{-1}$ onto $\widetilde{P}:=[0,\infty) \rtimes (0,1]$. Thus we have $\mathcal{W}_{\ell}(P) \simeq \mathcal{W}_{r}(P^{-1}) \simeq \mathcal{W}_{r}(\widetilde{P})$.

Let $\widetilde{X}:=[-\infty,0] \times [1,\infty] \subset Y$. Then $\widetilde{X}$ is $\widetilde{P}$ invariant. We leave it to the reader to check that $\widetilde{X}_{0}:=\widetilde{X}Int(\widetilde{P})=[-\infty,0)\times (1,\infty]$ which is open in $Y$. As before one can show that $\mathcal{W}_{r}(\widetilde{P})$ is isomorphic to $C_{red}^{*}(\widetilde{X} \rtimes \widetilde{P})$ and Morita equivalent to $C_{0}(Y) \rtimes (\mathbb{R} \rtimes (0,\infty))$. It would be interesting to decide whether $\mathcal{W}_{\ell}(P)$ and $\mathcal{W}_{r}(P)$ are isomorphic or not. \\

\nocite{Nica_WienerHopf}
\nocite{Renault_Muhly}
\nocite{Renault_thesis}
\nocite{Ore}
\nocite{MRW}
\bibliography{references}
\bibliographystyle{amsalpha}

\noindent{\sc Jean Renault} (\texttt{Jean.Renault@univ-orleans.fr})\\
         {\footnotesize Universit\'e d'Orl\'eans et CNRS (UMR 7349 et FR 2964), \\
D\'epartement de Math\'ematiques, 45067, Orl\'eans Cedex 2, France. }\\[1ex]
{\sc S. Sundar}
(\texttt{sundarsobers@gmail.com})\\
         {\footnotesize  Chennai Mathematical Institute, H1 Sipcot IT Park, \\
Siruseri, Padur, 603103, Tamilnadu, INDIA.}

\end{document}